\documentclass[11pt]{article}
\usepackage{amsfonts}
\usepackage{amssymb}
\usepackage{latexsym}
\usepackage{amsmath}
\usepackage{fullpage}
\usepackage[all]{xy}
\font\velkymat=cmmi10 scaled 4700
\font\velkymatt=cmmi10 scaled 3500

\font\strednimat=cmmi10 scaled 2500
\def\odstavec#1{\smallskip\noindent {\bf #1.}}

\newcommand\ulap[1]{\vbox\@to\z@{{\vss#1}}}%
\newcommand{\su}{\mathfrak {su}}
\newcommand{\Su}{SU}
\newcommand{\un}{\mathfrak {u}}
\newcommand{\Un}{U}
\renewcommand{\sp}{\mathfrak {sp}}
\newcommand{\g}{\mathfrak g}
\newcommand{\h}{\mathfrak h}
\newcommand{\I}{\mathbb I}

\newcommand{\Ad}{\operatorname{Ad}}
\newcommand{\ad}{\operatorname{ad}}
\newcommand{\hookarrow}{\hookrightarrow}
\newcommand{\diag}{\operatorname{diag}}
\newcommand{\D}{\mathfrak D} 
\newcommand{\C}{\mathbb C}
\newcommand{\R}{\mathbb R}
\newcommand{\Z}{\mathbb Z}
\newcommand{\End}{\operatorname{End}}
\newcommand{\w}{\omega}
\newcommand{\del}{\partial}

\newcommand{\Cof}{\operatorname{Cof}}
\newcommand{\T}{{\mbox{\rm T}}}
\newcommand{\ov}{\overline}
\newcommand{\wh}{\widehat}

\newcommand{\tr}{\operatorname{tr}}
\newcommand{\vloz}{\hookrightarrow}
\newcommand{\X}{\mathfrak X}
\newcommand{\bx}{{\bf x}} 
\newcommand{\by}{{\bf y}}
\newcommand{\on}{\overline\nabla}
\def\ol#1{\overline #1}
\newcommand{\hn}{\widehat \nabla}
\newcommand{\oX}{\overline X}
\newcommand{\hX}{\hat X}
\newcommand{\oY}{\overline Y}
\newcommand{\hY}{\hat Y}
\newcommand{\og}{\overline g}

\newcommand{\oZ}{\overline Z}
\newcommand{\hZ}{\hat Z}
\newcommand{\oR}{\overline R}
\newcommand{\spa}{\operatorname{span}}
\newcommand{\gr}{\operatorname{gr}}
\renewcommand{\Re}{\operatorname{Re}}
\def\cC{{\cal C}}
\def\odst #1{\section{#1}}
\newenvironment{dk}{\medskip \noindent{\bf Proof.}}{\hbox{}\hfill
$\Box$\bigskip}

\newtheorem{thm}{Theorem}[section]
\newtheorem{prop}[thm]{Proposition}
\newtheorem{lem}[thm]{Lemma}

\newtheorem{rem}[thm]{Remark}
\newtheorem{ex}[thm]{Example}

\def\ov{\overline}
\def\ot{\otimes}

\def\ra{\rightarrow}

\newcommand{\U}{{\mbox{\rm U}}}

\newcommand{\G}{{\mbox{\rm G}}}

\newcommand{\eH}{{\mbox{\rm H}}}
\newcommand{\eP}{{\mbox{\rm P}}}

\newcommand{\impl}{\Longrightarrow}
\newcommand{\al}{\alpha}
\newcommand{\be}{\beta}

\newcommand{\la}{\lambda}
\newcommand{\om}{\omega}
\newcommand{\Om}{\Omega}
\renewcommand{\th}{\theta}

\renewcommand{\phi}{\varphi}

\def\Z{{\Bbb Z}}
\def\R{{\Bbb R}}
\def\C{{\Bbb C}}

\def\P{{\Bbb P}}

\def\F{{\Bbb F}}

\def\V{{\cal V}}

\def\cC{{\cal C}}
\def\cR{{\cal R}}
\def\hR{\widehat R}

\def\su{{\frak {su}}}

\def\sl{{\frak {sl}}}
\def\sp{{\frak {sp}}}

\def\g{{\frak g}}
\def\h{{\frak h}}

\def\p{{\frak p}}

\renewcommand{\l}{{\frak l}}
\def\X{{\frak X}}

\renewcommand{\frak}{\mathfrak}
\renewcommand{\Bbb}{\mathbb}
\def\be{\begin{equation}}
\def\ee{\end{equation}}
\def\bi{\begin{enumerate}}
\def\ei{\end{enumerate}}
\def\ba{\begin{array}}
\def\ea{\end{array}}
\def\bea{\begin{eqnarray}}
\def\eea{\end{eqnarray}}
\def\ben{\begin{enumerate}}
\def\een{\end{enumerate}}

\def\hook{\mbox{}\begin{picture}(10,10)\put(1,0){\line(1,0){7}}
 \put(8,0){\line(0,1){7}}\end{picture}\mbox{}}

\title{{\bf Bochner-Kaehler metrics and connections of Ricci type}}

\def\ftnote#1{\def\footnotemark{}\footnote{#1}\setcounter{footnote}{0}}

\author{Martin Pan\'ak${}^{1,3}$\ftnote{\kern-4.5pt${}^{1}$
Institute of Mathematics, Academy of Sciences of the Czech Republic. e-mail: 
naca@math.muni.cz}
\mbox{\hspace{1cm}} Lorenz
J. Schwachh\"ofer${}^{2,3}$
\ftnote{\kern-4.5pt${}^{2}$Universit\"at 
Dortmund, Vogelpothsweg 87, 44221 Dortmund, Germany. e-mail: 
lschwach@math.uni-dortmund.de}
\ftnote{\kern-4.5pt${}^{3}$Both authors were supported by the Schwer\-punkt\-pro\-gramm Globale 
Differentialgeometrie of the Deutsche Forschungsgesellschaft, the first one also by the
grant nr. 201/05/P088 of the Grant academy of the Czech republic}}

\begin{document}
\maketitle

\begin{abstract}
We apply the results from \cite{CS} about special symplectic geometries to the case of Bochner-Kaehler metrics.
We obtain a (local) classification of these based on the orbit types of the adjoint action in $\su(n,1)$. 
The relation between Sasaki and Bochner-Kaehler metrics in cone and transveral metrics constructions 
is discussed. The connection of the special symplectic and Weyl connections is outlined. The duality between
the Ricci-type and Bochner-Kaehler metrics is shown.
\end{abstract}


\odst{Bochner-Kaehler metrics} The curvature tensor of the
Levi-Civita connection of
a Kaehler metric $g$ decomposes (under the action of $\un(n)$) into its Ricci and Bochner part
(\cite{Bo}).
The metric is said to be Bochner-Kaehler, 
iff the Bochner part  of its curvature tensor vanishes.

A remarkable relationship was revealed among following types of
geometric structures in the article \cite{CS}: manifolds with a connection of Ricci type, manifolds
with a connection with the special symplectic holonomy,  pseudo-Riemannian
Bochner-K\"ahler
structures, manifolds with a Bochner-bi-Lagrangian connection.
  All these geometric objects are instantons of the same construction, and
they are called special symplectic geometries.
The word "symplectic" comes from the fact that they all carry a symplectic connection; special stands
for the common special  type
of the curvature of the connection: let
$(M,\w)$ be a symplectic manifold.
Then the curvature of the special
symplectic geometries is of the form
\begin{eqnarray}
\label{krivost}
R_{h}(X,Y)=2\w(X,Y)h+X\circ(hY)-Y\circ(hX),
\end{eqnarray}
where $\h\subset\sp(n,\R)$ (or $\sp(n,\C)$) is a Lie algebra, $h\in\h$,
$\circ:S^2(TM)\to \h$, is an $\h$-equivariant product
with
special properties (see \cite{CS}). 
For
Bochner-Kaehler structures the special form of curvature translates as follows:
let $(M,g,J,\w)$ be a Kaehler structure on a manifold $M$. That is $J$ is
the orthogonal complex structure which is parallel with respect to the
Levi-Civita connection of $g$, and the Kaehler form $\w$ is defined by
$\w(x,y)=g(x,Jy)$. The Kaehler structure  is Bochner-Kaehler iff the curvature of the Levi-civita connection
of the metric has the above form, where $\h=\un(n)$, and $\circ$ is given as:

\begin{eqnarray} 
\label{Bochnerkolecko}
(X\circ Y)Z=\w(X,Z)Y+\w(Y,Z)X+\w(JX,Z)JY+\w(JY,Z)JX+\w(JX,Y)JZ.
\end{eqnarray}
That is iff the curvature is of the form
\begin{eqnarray}
\label{krivost2}
  R_\rho(X,Y)&=&2g(X,JY)\rho+2g(X,\rho Y)J+(\rho Y\wedge JX)\nonumber\\
& &-(\rho X\wedge
   JY)+(X\wedge J\rho Y)-(Y\wedge J\rho X),
\end{eqnarray}
where $(X\wedge Y)Z=g(X,Z)Y-g(Y,Z)X$.

\odst{Kaehler and Sasaki manifolds} The following considerations are motivated by the
lecture given by  Krzysztof Galicki at Winter School of Geometry and
Physics, Srn\'\i, Czech republic, 2004. (\cite{BG}).

\odstavec{Sasaki metric} One of the possible (equivalent) definition of the {\it Sasakian  manifold} $C$ 
is that it is a Riemannian manifold with the metric $g$, on which there exists 
a unit length Killing vector field $\xi$ such that the curvature tensor $R$ of the
Levi-Civita connection $\nabla$ of $g$ satisfies:
$$
R(X,\xi)Y=g(\xi,Y)X-g(X,Y)\xi.
$$
The one form  $\lambda$ dual to $\xi$  defines a contact distribution
$\D=\{X\in TC| \lambda(X)=g(X,\xi)=0\}$. The vector field $\xi$ is called
the characteristic vector field of the contact distribution $\D$.

\odstavec{Transversal Kaehler metric}
Further consider  $J$ 
defined by
$J(X)=-\nabla_X\xi$. It is an automorphism of the tangent bundle $TC$
and its restriction to $\D$ gives rise to a 
complex structure $J$ on $\D$. Then $(\nabla_XJ)(Y)=0$ for $X$, $Y\in \D$ and
thus there is a so called {\it transversal Kaehler structure} on $\D$. The Kaehler structure
then factorizes to the set of leaves of the foliation generated by $\xi$
(the characteristic foliation)
if this is locally an orbifold. See \cite{BG} for details.

Conversely, the tranversal Kaehler structure on a compact distribution $\D$
(given as $\D=\{X|\lambda(X)=0, \lambda\in
\Omega^1(C)\}$)
on a manifold $C$ translates to a Sasakian structure on $C$:
given a metric $g_{\D}$ on $\D$ with a parallel complex structure $J$ on $\D$,
and a transversal symmetry $\xi$ of $\D$, one
extends $g_{\D}$ to the whole of $TC$ with
$g(X,Y)=g_{\D}(X,Y)$ for $X$, $Y\in \D$, $g(\xi,\xi)=1$ and $g(\xi,X)=0$.

Unlike a general Kaehler metric, any Bochner-Kaehler metric can be realized (locally)
as the transverse metric of an appropriate Sasaki metric.
The Theorem B from \cite{CS} (see also \ref{locallygenerating}) says, that given a simply connected manifold $M$ with a
Bochner-Kaehler metric $g$ there is always
a principal $T$-bundle $\pi: C\to M$, where $T$ is a one-dimensional Lie
group, and this bundle carries a connection whose curvature equals $-2\w$,
where $\w$ is the Kaehler form corresponding to $g$.
The horizontal distribution of the connection  yields a
contact distribution on $C$. Thus we have following:

\begin{prop}
\label{existenceofSasaki}
Let $M$ be a $2n$-dimensional (real) manifold with Bochner-Kaehler metric $g$, 
$J$ be the corresponding complex structure.
Then there exists a Sasaki manifold such that the set of leaves of the
characteristic foliation is isomorphic (together with from the Sasaki one induced
structure) to some cover of $M$ (with the Bochner-Kaehler structure induced from $M$).  
\hbox{}\hfill$\Box$
\end{prop}

\odstavec{Cone metric} On the other hand, a manifold $C$ is Sasakian if and only if 
the "cone metric" $(t^2\cdot g+ (dt)^2)$ on $C\times \R_+$ is Kaehler,
where the complex structure $J'$ on the cone is the extension of $J$ such that 
\begin{equation}
\label{Jstruktura}
J'(\xi)=t\del_t, \quad J'(\del_t)=-\frac{1}{t}\xi,
\end{equation}
$\xi\in \X(\widehat C)$ being the lift of the characteristic vector field $\xi$ on $C$.

Following the ideas from \cite{BG}, there arise  questions 
what happens if we require the transversal metric to be Bochner-Kaehler.
What special has to be the Sasaki metric on $C$?
Will then the "cone metric" be Bochner-Kaehler?

Let us state a technical lemma about Bochner-Kaehler manifolds.
\begin{lem}
\label{directionflat}
Let $M$ be a manifold with a Bochner-Kaehler metric $g$ and let $X_0$ be a 
non-vanishing vector
field on $M$ such that  the curvature $R$ of the Levi-Civita connection
satisfies $R(X_0,JX_0)\equiv 0$. Then $M$ is flat.
\end{lem}
\begin{dk}
The curvature of the Levi-Civita connection of the Bochner-Kaehler metric is
of the form $R_\rho$, $\rho\in\un(n)$ (see $\eqref{krivost2}$), the vector
field $X_0$ can be normed to the unit length and we can write
\begin{eqnarray*}
0=R(X_0,JX_0)=-2\rho+2g(X_0,\rho JX_0)J-2X_0\wedge \rho X_0-2JX_0\wedge\rho JX_0,
\end{eqnarray*}
that is
\begin{eqnarray}
\label{ro}
\rho=g(X_0,\rho JX_0)J-X_0\wedge\rho X_0-JX_0\wedge\rho JX_0.
\end{eqnarray}
Applying $\rho$ to $X_0$ we get
\begin{eqnarray*}
\rho X_0&=&g(X_0,\rho JX_0)JX_0-g(X_0,X_0)\rho X_0+g(\rho X_0,X_0)X_0-g(JX_0,X_0)\rho JX_0+g(\rho JX_0,X_0)JX_0\\
&=&2g(X_0,\rho JX_0)JX_0-\rho X_0, 
\end{eqnarray*}
and we have
\begin{equation}
\label{ro2}
\rho X_0=g(X_0,\rho JX_0)JX_0.
\end{equation}
that means $\rho X_0=c JX_0$ for a real valued
function $c$ on $M$. Substituing
back to $\eqref{ro2}$ we get
$$
cJX_0=g(X_0,J(cJX_0))J=-cJX_0,
$$
and $c=0$. That is $\rho X_0=0$, and $0=\rho JX_0=J\rho X_0$, and
the formula $\eqref{ro}$ implies $\rho=0$, that is the curvature vanishes.
\end{dk}

\begin{prop} 
\label{Cone construction}
Let $M$ be a $2n$-dimensional (real) manifold with Bochner-Kaehler metric $g$, 
$J$ be the corresponding complex structure, further let $C$
be the Sasakian manifold from the theorem \eqref{existenceofSasaki} with the dimension $2n+1$.
Then the manifold
$\widehat C= C\times \R_+$ with the complex structure defined by
$\eqref{Jstruktura}$ 
is Bochner Kaehler if and only if $M$ is locally isomorphic (as the Kaehler structure)
to the complex projective space $\C P^n$. The cone $\widehat C$ is then a flat manifold. 
\end{prop}

\begin{dk}
The existence of the given manifold is just consequence of the proposition \cite{existenceofSasaki} and the
cone construction.
We will write $\overline g$ and $\hat g$ for the metrics on $C$ and $\widehat C$ respectively, and
$\on$, $\hn$ for the corresponding Levi-Civita connection. The curvature of $\kappa$ gives us
$$[\oX,\oY]=\overline[X,Y]-2\w(X,Y)\xi,$$
for $X$, $Y\in\X(M)$ and since $[\xi,X]=0$ there is $\on_\xi X=\on_X\xi$.
That gives for the 
torsion of the connection $\on$
\begin{eqnarray*}
0=T(\oX,\oY)&=&\on_{\oX}\oY-\on_{\oY}\oX-[\oX,\oY]=\on_{\oX}\oY-\on_{\oY}\oX-\overline{[X,Y]}+2\w(X,Y)\xi\\
&=&\on_{\oX}\oY-\on_{\oY}\oX-\overline{\nabla_XY}+\overline{\nabla_YX}+2\w(X,Y)\xi
\end{eqnarray*}
On the other hand
\begin{eqnarray*}
0=\xi\og(\oX,\oY)&=&\og(\on_{\xi}\oX,\oY)+\og(\oX,\on_\xi\oY)\\
&=&X(\og(\xi,\oY))-\og(\xi,\on_{\oX}\oY)+\oY(\og(\oX,\xi))
-\og(\on_{\oY}\oX,\xi)\\
&=&-\og(\xi,\on_{\oX}\oY+\on_{\oY}\oX).
\end{eqnarray*}
We conclude that $\on_{\oX}\oY+\on_{\oY}\oX=0$ and
\begin{equation}
\on_{\oX}\oY=\overline{\nabla_XY}-\w(X,Y)\xi.
\end{equation}
Further
\begin{eqnarray*}
0=\oX \og(\xi,\xi)=2\og(\on_{\oX}\xi,\xi),
\end{eqnarray*}
and
\begin{eqnarray*}
0=\oX\og(\xi,\oY)&=&\og(\on_{\oX}\xi,\oY)+\og(\xi,\on_{\oX}\oY)\\
                 &=&\og(\on_{\oX}\xi,\oY)-\w(X,Y)\\
                 &=&g(\pi(\on_{\oX}\xi),Y)+g(JX,Y).\\
\end{eqnarray*}
That is
\begin{equation}
\on_\xi\oX=\on_{\oX}\xi=\overline{-JX}.
\end{equation}
With the like-wise computation one gets $\on_{\xi}\xi=0$.

\noindent
Similarly one gets for $X$, $Y\in \X(C)$
\begin{equation*}
\hn_{\hX}\hY=\widehat{\on_XY}-tg(X,Y)\del_t,\ \ \mbox{especially}\ \ \hn_{\hat{\xi}}\hat\xi=-t\del_t,
\end{equation*}
and

\begin{equation*}
\hn_{\hX}\del_t=\hn_{\del_t}\hX={1\over t}\hX, \ \mbox{especially}\hn_{\del_t}\hat\xi={1\over t}\hat\xi\ \ \mbox{and}\ \ \hn_{\del_t}\del_t=0.
\end{equation*}

There is the canonical projection $\pi:\widehat C=C\times\R_+\to C \to M=C/T$. 
From the above two equations we see, that the fibres of this projection sit totally geodesicly in $\widehat C$ and
the metric on the fiber gives the "true" 2-dimensional cone. That is the restriction of the metric to the fibers
is flat and we have $\hR(\del_t,J(\del_t)){\hat\xi}=0$ and $\hR(\del_t,J(\del_t))\del_t$.
It remains to compute $\hR(\del_t,J(\del_t)$ on the lifts of the vectors in $TM$.
We will
write $\widehat X$ for a lift of a vector $X\in TM$ with respect to the projection $\pi$.
\begin{eqnarray*}
\hR(\del_t,J(\del_t))\hat X&=&-{1\over t}\left(\hn_{\del_t}\hn_{\hat\xi}\hX-\hn_{\hat\xi}\hn_{\del_t}\hX
-\hn_{[\del_t,\hat\xi]}\hX\right)\\
 &=&-{1\over t}\left(\hn_{\del_t}(-\widehat{JX})-\hn_{\hat\xi}({1\over t}\hX)\right)\\
 &=&-{1\over t}\left(-{1\over t}\widehat{JX}+{1\over t}\widehat{JX}\right)=0.
\end{eqnarray*}
Consequently $\hR(\del_t,J(\del_t))=0$ and according to the previous lemma \ref{directionflat} the metric on $\widehat C$ is flat.

There is the following known relation 
between the curvature $\overline R$
of  the metric $\og$ on $C$ and $\hat R$, the curvature of the Levi-Civita connection $\hn$
of the cone metric on $\widehat C$:
\begin{eqnarray*}
\hat R(\hX,\hY)\hZ=\overline R(X,Y)Z+\og(X,Z)Y-\og(Y,Z)X,
\end{eqnarray*}
for $X$,$Y$, $Z\in T(C)$, which can be checked with an easy computation excercise.

This shows, that $C$ is a manifold with constant sectional curvature $K=1$
(see \cite{KN}), that is locally a unit sphere. 

Now the vector field $\xi$ on $C$ is from the construction a Killing one,
and has the constant unit length
that is the leaves of the foliation are the circles coming
from the natural $\C$ action on the sphere (at least locally; the unit 
length of the Killing field excludes other circle actions). Then the
resulting factor space, that is $M$, is locally isomorphic to the complex
projective space.  
\end{dk}

\odst{General construction}

Let us quickly review the construction from \cite{CS}, which gives rise to all special symplectic geometries.
All manifolds with special symplectic connection are locally isomophic
to the factor manifold of the oriented projectivization of the cone $\cC=Ad_Gx\subset \g$, where $x$ is an apropriate element
in the parabolic 2-gradable Lie algebra $\g$, where we factor along the flow of the convenient vector field.
The special symplectic connection is then induced on the factor from one of the components of the Maurer-Cartan form 
on $\g$, which decomposes due to the 2-grading. 

Some of the standard notions
from the theory of contact structures are used without definitions.
The reader can find them and all the proofs of the theorems stated in this section, in \cite{CS}.

\odstavec{Symplectic algebra as subalgebra of a 2-graded algebra}
Let $V$ be a vector space (either real or complex) with a symplectic form $\w$.
Let $\h\subset\sp(V,\w)=\{h\in \End(V)|\w(x,y)+\w(x,hy)=0 \mbox{ for all }x,y \in V\}$ such, that there exists an 
$\h$-equivariant map $\circ: S^2(V)\to \h$ and an $\ad_\h$-invariant inner product $(,)$ which satisfy
the following identities:
\begin{eqnarray*}
(h, x \circ y) &=& \om(hx, y) = \om(hy, x)\\
(x \circ y) z - (x \circ z) y &=& 2\ \om(y,z) x - \om(x,y) z + \om(x,z) y,
\end{eqnarray*}
for all $x,y,z \in V$ and $h \in \h$.

Then there exists a unique simple Lie algebra $\g$ with a 2-grading of the parabolic type, that is
$$\g=\g^{-2}\oplus \g^{-1}\oplus \g^0\oplus\g^1\oplus\g^2,$$
where $\g^{-2}$ and $\g^2$ are one-dimensional. 
The grading corresponds to $\h$ in the following sense:
$$
\g^{ev} := \g^{-2} \oplus \g^0 \oplus \g^2 \cong \sl_{\al_0} \oplus 
\h\ \ \ \ \mbox{and}\ \ \ \ \g^{odd} := \g^{-1} \oplus \g^1 
\cong \F^2 \ot V\ \ \ \ \mbox{as a $\g^{ev}$-module},
$$
where $\g^{2}$, resp. $\g^{-2}$, are root spaces of a long root $\al_0$, resp $-\al_0$,
and $\sl_{\al_0}$ is the Lie algebra isomorphic to $\sl(2,\F)$ generated by the
root spaces and the corresponding coroot $H_{\al_0}$ which lies in $\g^0$.
We will also write $\p=\g^0\oplus\g^1\oplus\g^2$ for the parabolic subalgebra of $\g$
and $\p_0:=\h\oplus\g^1\oplus\g^2$. Let further $P$ and $P_0$ be corresponding
connected subgroups of $G$.

Further we fix a non-zero $\F$-bilinear area form $a \in \Lambda^2 (\F^2)^*$. There 
is a canonical $\sl(2,\F)$-equivariant isomorphism
\be \label{eq:define sl(2)}
S^2(\F^2) \longrightarrow \sl(2,\F),\ \ \ \ (ef) \cdot g := a (e, g) f + a(f, 
g) e\ \ \mbox{for all $e,f,g \in \F^2$},
\ee
and under this isomorphism, the Lie bracket on $\sl(2,\F)$ is given by
\be \label{eq:bracketsl2}
{}[ef, gh] = a(e, g) fh + a(e, h) fg + a(f, g) eh + a(f, h) eg.
\ee
Thus, if we fix a basis $e_+, e_- \in \F^2$ with $a(e_+, e_-) = 1$, then we 
have the identifications
\[
H_{\al_0} = - e_+ e_-,\ \ \ \g^{\pm2} = \F e_\pm^2,\ \ \ \g^{\pm 1} = e_\pm 
\ot V.
\]

\odstavec{The cone in 2-gradable algebra and its projectivization}
Using the Cartan-Killing form (up to the multiple) we identify $\g$ 
and $\g^*$, and we define the root cone $\hat \cC$ and its 
(oriented) projectivization $\cC$ as follows:
\[
\hat \cC := \G \cdot e_+^2 \subset \g \cong \g^*,\ \ \ \ \ \ 
\cC := p(\hat \cC) \subset \P^o(\g) \cong \P^o(\g^*),
\]
where $\P^o(\g)$ is the set of {\em oriented} lines in $\g$, i.e. $\P^o 
\cong S^d$ if $\F = \R$, and $\P^o \cong \C\P^d$ if $\F = \C$, where $d = 
\dim \g - 1$, and where $p: \g \backslash 0 \ra \P^o(\g)$ is the principal 
$\R^+$-bundle ($\C^*$-bundle, respectively) defined by the canonical 
projection. Thus, the restriction $p: \hat \cC \ra \cC$ is a principal
bundle as well.

\odstavec{Contact structure on the projectivized cone}
Being a coadjoint orbit, $\hat \cC$ carries a canonical $\G$-invariant 
symplectic structure $\Om$. Moreover, the {\em Euler vector field} 
defined by
\[
E_0 \in {\frak X}(\hat \cC),\ \ \ \ \ (E_0)_{v} := v
\]
generates the principal action of $p$ and satisfies ${\frak L}_{E_0}(\Om) = 
\Om$, so that the distribution 
\begin{equation}
\label{Kontaktnidistribuce}
{{\cal D}} = dp(E_0^{\perp_\Om}) \subset T\cC
\end{equation}
yields a $\G$-invariant contact distribution on $\cC$, see \cite{CS}, Proposition 3.2.

\odstavec{The cone as homogeneous space}
Let $\lambda=\iota_{E_0}(\Omega)$. Then we define the bundle $\cR$:
\[
{\frak R} := \{ (\la, \hat \xi) \in \hat \cC \times T\hat \cC \subset T^*\cC 
\times T\hat \cC \mid \la(dp(\hat \xi)) = 1\}.
\]
Let $\eP$ and $\eP_0$ be subgroups of $G$ corresponding to the subalgebras $\p$ and $\p_0$ of $\g$.

\begin{lem} 
\label{Homogeneous spaces}
As homogeneous spaces, we have $\cC = \G/\eP$, $\hat \cC = 
\G/\eP_0$, and ${\frak R} = \G/\eH$.
\end{lem}
\begin{dk}
See \cite{CS}
\end{dk}

{\bf Transversal symmetry defines the geometry}

For each $a \in \g$ we define the vector fields $a^* \in {\frak 
X}(\cC)$ and $\hat a^* \in {\frak X}(\hat \cC)$
corresponding to the infinitesimal action of $a$, i.e.
\be \label{eq:xi0}
(a^*)_{[v]} := \left. \frac d{dt} \right|_{t=0} (\exp(t a) \cdot [v])
\ \ \ \ \ \mbox{and}\ \ \ \ 
(\hat a^*)_v := \left. \frac d{dt} \right|_{t=0} (\exp(t a) \cdot v).
\ee
Note that $a^*$ is a contact symmetry (with respect to the canonical contact distribution on $\hat \cC$),
and $\hat a^*$ is its 
Hamiltonian lift. Let
\be \label{eq:cC0}
\hat \cC_a := \{ \la \in \hat \cC \mid \la(a^*) \in \R^+ (\in 
\C^*, \mbox{ respectively})\}\ \ \ \ \ 
\mbox{and}\ \ \ \ \ \cC_a := p(\hat \cC_a) \subset \cC,
\ee
so that $p: \hat \cC_a \ra \cC_a$ is a principal $\R^+$-bundle 
($\C^*$-bundle, respectively) and the restriction of $a^*$ to $\cC_a$ is 
a positively transversal contact symmetry. Then there exists a unique section $\la$ of the bundle
$p: \hat \cC_a \ra \cC_a$ such that $\la(a^*)=1$ and 
therefore, we obtain the section 
\be \label{eq:sectionR}
\sigma_a: \cC_a \longrightarrow {\frak R}=\G/\eH,\ \ \ \ \ \sigma_a(u) := (\la(u), \hat a^*(u)) \in {\frak R}.
\ee

Let $\pi: \G \ra \G/\eH = {\frak R}$ be the canonical projection, and let 
$\Gamma_a:= \pi^{-1}(\sigma_a(\cC_a)) \subset \G$. The 
restriction $\pi: \Gamma_a \ra \sigma_a(\cC_a) \cong \cC_a$ is then a 
principal $\eH$-bundle.

\begin{thm} \label{thm:Maurer-Cartan}
Let $a \in \g$ be such that $\cC_a \subset \cC$ from (\ref{eq:cC0}) is 
non-empty, define $a^* \in {\frak X}(\cC)$ and $\hat a^* \in 
{\frak X}(\hat \cC)$ as in (\ref{eq:xi0}), and let $\pi: \Gamma_a \ra \cC_a$ 
with $\Gamma_a \subset \G$ be the principal $\eH$-bundle from above. Then 
there are functions $\rho : \Gamma_a \ra \h$, $u : \Gamma_a \ra V$, $f: 
\Gamma_a \ra \F$ such that 
\be \label{eq:formA}
\Ad_{g^{-1}}(a) = \frac12 e_-^2 + \rho + e_+ \ot u + \frac12 f e_+^2
\ee 
for all $g \in \Gamma_a$. 
\end{thm}

The restriction of the $\mu_\h + \mu_{-1} + \mu_{-2}$ part of the Maurer-Cartan form 
to $\Gamma_a$ yields a pointwise linear isomorphism 
$T\Gamma_a \ra \h \oplus \g^{-1} \oplus \g^{-2}$, and we can further decompose it as

\[
\mu_\h + \mu_{-1} + \mu_{-2} = -2 \kappa\ \left(\frac12 e_-^2 + 
\rho\right) + e_- \ot \th + \eta,\ \ \ \ 
\kappa \in \Om^1(\Gamma_a),\ \ \th \in \Om^1(\Gamma_a) \ot V,\ \ \ 
\eta \in \Om^1(\Gamma_a) \ot \h.
\]

\begin{thm} \label{thm:canonicalconn}
Let $a \in \g$ and $\cC_a \subset \cC$ as before. Let $U \subset \cC_a$ be a
regular open subset , i.e. the local quotient $M_U := \T_a^{loc} \backslash 
U$ is a manifold, where
\[
\T_a := \exp(\F a) \subset \G.
\]
Let $\om \in \Om^2(M)$ be the unique symplectic form on $M_U$, such that $\pi^*(\om)=-2d(E_0\hook \Om)$. Then 
$M_U$ carries a canonical special symplectic connection associated to $\g$,
and the (local) principal $\T_a$-bundle $\pi: U \ra M$ admits a connection 
$\kappa \in \Om^1(U)$ whose curvature is given by $d\kappa = \pi^*(\om)$.
\end{thm}

\begin{dk}
Sketch: The connections forms of the desired connections are projections of the forms
$\eta$ and $\kappa$ on $\Gamma_a$ over U to the corresponding $H$-bundle over $M_U$.  
\end{dk}

Conversely, any manifold with special symplectic connection comes in this way (locally). Namely, there is the following
theorem (Theorem B from \cite{CS}):

\begin{thm}
\label{locallygenerating}
Let $(M,\om)$ be a symplectic manifold with a special symplectic connection of class $C^4$, and let $\g$
be the Lie algebra associated to the special symplectic condition as described at the beginning of
this section.
\begin{itemize}
\item[i)] Then there is a principal $\hat T$-bundle $\hat M\to M$, where $\hat T$ is a one dimensional Lie group
which is not necesarily connected, and this bundle carries a principal connection with curvature $\omega$.
\item[ii)] Let $T\subset\hat T$ be the identity component. Then there is an $a\in \g$ such that 
$T\equiv T_a\subset G$, and a $T_a$-equivariant local diffeomorphism $\hat \iota:\hat M\to C_a$ which for
each sufficiently small open subset $V\subset \hat M$ induces a connection preserving diffeomorphism $\iota:
T^{\mbox{\scriptsize loc}}\V\to T^{\mbox{\scriptsize loc}}\U=M_U$, where $U:=\hat\iota(V)\subset C_a$ and
$M_U$ carries the connection from \ref{thm:canonicalconn}.
\end{itemize}
\hbox{}\hfill
$\Box$\bigskip
\end{thm}

\odst{Construction of Bochner-Kaehler metrics}
\subsection{A little of linear algebra}

Let us first recall some facts from the linear algebra. Let $V$ be a complex $(n+1)$-dimensional
vector space, $h$ a hermitian form of the signature $(n,1)$ on $V$. The (real) Lie algebra $\un(n,1)$ is
defined as follows:

$$\un(n,1):=\{A\in\End(V)|h(Av,w)+h(v,Aw)=0, v,w\in V\}.$$
It is the Lie alebra of the Lie group
$$\Un(n,1)=\{A\in Aut(V)|h(Av,Aw)=h(v,w), v,w\in V\}.$$

Take matrices with the determinant 1 in $\Un(n,1)$ or traceless matrices in $\un(n,1)$ to get something special.

\begin{ex}
\label{Gradace algebry}
We consider the standard hermitian form on the complex space $\C^{n+1}$ of signature $(n,1)$
($h(\bx,\by)=\sum_{i=1}^n x_i\overline y_i - x_{n+1}\overline y_{n+1}$, for $\bx$, $\by\in C^{n+1}$. Then as a matrix 
algebra $\su(n,1)$ can be written as
\begin{equation*}
\su(n,1)=\left\{\left(\begin{array}{c|c}
\hphantom{aa}\hbox{\LARGE A}\hphantom{aa}&\raise .1cm\hbox{\bf v}\\
\hline
\bf v^*& -\tr A\\
\end{array}\right), A\in\un(n), 
{\bf v}=\left(\begin{array}{c}v_1\\ 
\vdots\\ 
v_n
\end{array}\right)\in \C^n\right\}.\\
\end{equation*}
The elements of the bundle $\Gamma_a$ from the general construction in the previous section,
can be described with the structure functions from \eqref{thm:Maurer-Cartan}
as follows:
\begin{equation}
\label{matrixform}
\begin{pmatrix}
\rho-\frac{1}{n+2}(\tr\rho )\I_n&u&u\\
-u^*&-\frac12(\tr\rho-i(f+1))&\frac i2(f-1)\\
u^*&\frac i2(1-f)& -\frac12(\tr\rho+i(f+1))
\end{pmatrix}.
\end{equation}
The grading of $\su(n,1)=\g^{-2}\oplus \g^{-1}\oplus \g^0\oplus \g^1\oplus \g^2\equiv 
\R e_-^2\oplus e_-\otimes V\oplus (\h\oplus \R e_+ e_-)\oplus e_+\otimes V\oplus \R e_+^2$ 
 is given as follows (all the matrices are $(n+1)\times (n+1)$ ones):
\begin{eqnarray*}
e_{\pm}^2&=&
\left(\begin{array}{c|cc}
\vbox to 0pt{\hbox{\velkymat 0}}&0&0\\
&\vdots&\vdots\\
&0&0\\
\hline
0\ldots 0&i& \pm i\\
0\ldots 0&\mp i&-i\\
\end{array}\right),\qquad
H_{\al_0}=e_+e_-=
\left(\begin{array}{c|cc}
\vbox to 0pt{\hbox{\velkymat 0}}&0&0\\
&\vdots&\vdots\\
&0&0\\
\hline
0\ldots 0&0&1\\
0\ldots 0&1&0\\
\end{array}\right),\\
\smallskip
\h&=&\scriptsize\left\{\left(\begin{array}{c|rc}
\vbox to 0pt{\hbox{\velkymatt A}}&0&0\\
&\vdots&\vdots\\
&0&0\\
\hline
0\ldots0&\rlap{\tiny-$\frac12\tr A$}\hphantom{000}& 0\\
0\ldots 0&0&\hphantom{0000}\llap{\tiny-$\frac12\tr A$}\\
\end{array}\right), A\in\un(n-1)\right\},\\
\smallskip
\g^{\pm 1}&=&\left\{
\left(\begin{array}{c|cc}
\vbox to 0pt{\hbox{\velkymat 0}}&&\\
&\bf v&\bf \pm v\\
&&\\
\hline
-\bf v^*&0& 1\\
\pm\bf v^*&1&0\\
\end{array}\right), 
{\bf v}=\left(\begin{array}{c}v_1\\ 
\vdots\\ 
v_{n-1}
\end{array}\right)\in \C^{n-1}\right\}.\\
\end{eqnarray*}
The action of the algebra $\h=\un(n-1)$ on $V$ is then given as the adjoint matrix action and one
easily computes that for $\rho\in \un(n-1)$ there is 
\begin{equation}
\label{akceun}
\rho\cdot u=\rho u + \frac12 \tr(\rho) u
\end{equation}

\end{ex}

The hermitian form $h$ is uniquely determined either by its real part $g$ (real valued symmetric bilinear form on $V$)
or by its imaginary part
$\omega$, the antisymmetric real valued form on $V$.
($\omega(x,y)=g(x,Jy)$, where $J$ is the complex structure on $V$).
\begin{lem}
\label{invariance}
There is a ${\Un(n,1)}$-equivariant  map
$m: V \to \un(n,1):\ x\mapsto x\wedge Jx$, 
where $$(x\wedge Jx)z=g(x,z)Jx-g(Jx,z)x,$$ 
and the $\Un(n,1)$ action on $\un(n,1)$ is given by $\ad$ representation.
\end{lem}
\begin{dk}
The morphism $(x\wedge Jx)$ is in $\un(n,1)$:

\begin{eqnarray*}
g((x\wedge Jx)y,z)+g(y,(x\wedge Jx)z)&=& g(g(x,y)Jx-g(Jx,y)x,z)+g(y,g(x,z)Jx-g(Jx,z)x)\\
&=&g(x,y)g(Jx,z)-g(Jx,y)g(x,z)+g(y,Jx)g(x,z)-\\
&&-g(y,x)g(Jx,z)\\
&=&0,
\end{eqnarray*}
for any $x$, $y$, $z\in V$. For the $A\in\Un(n,1)\subset GL(n+1,\C)$ there is:
\begin{eqnarray*}
\left(\Ad_A\circ m(x)\right)y&=&A\circ m(x)\circ A^{-1}(y)=A\circ\left(g(x,A^{-1})Jx-g(Jx,A^{-1}y)x\right)\\
&=&g(x,A^{-1}y)AJx-g(Jx,A^{-1}y)Ax=g(Ax,y)AJx-g(AJx,y)Ax\\
&=&g(Ax,y)JAx-g(JAx,y)Ax\\
&=&\left(m\circ \Ad_A(x)\right)(y),
\end{eqnarray*}
where we have used the invariance of $g$ with respect to the morphisms from $\Un(n,1)$.
\end{dk}

\begin{rem}
\begin{itemize}
\item[i)] The image of the morphism $m$ are from the definition rank one morphisms.
\item[ii)]
The value of the morphism $x\wedge Jx$ on a vector $z$ is actually $i\overline{h(x,z)}x$, but we stick
to write it in the form $g(x,z)Jx-g(Jx,z)x$, which comes from the morphism $x\wedge y$:
$(x\wedge y)z=g(x,z)y-g(y,z)x$.
\end{itemize}
\end{rem}

\begin{lem}
\label{Kruhova akce}
The morphism $m$ is  not injective: $x\wedge Jx$ and 
$y\wedge Jy$ determine the same element in $\un(n,1)$ if and only if
$x=e^{ik}y$, $k\in \R$, $x$, $y\in V$.
\end{lem}
\begin{dk}
If $x\wedge Jx$ and $y\wedge Jy$ determine the same morphism of $V$ then $x$ and $y$ 
lie on the same complex line, that is $y=ax+bJx$, $a$, $b\in \R$. Then there is
\begin{equation}
y\wedge Jy=(ax+bJx)\wedge(-bx+aJx)=(a^2+b^2)x\wedge Jx,
\end{equation} 
that is $x$ and $y$ differ by a multiple of a complex unit.
\end{dk}

\begin{lem}
\label{Bezestopy}
A morphism $x\wedge Jx$ lies in $\su(n,1)$ iff $g(x,x)=0$. 
\end{lem}
\begin{dk}
The trace of the rank one morphism is equal to the eigenvalue of a non-zero eigenvector which lies
in the image line. For the morphism $x\wedge Jx$ we take the eigenvector $x$. Then
\begin{equation}
(x\wedge Jx)x=g(x,x)Jx-g(Jx,x)x=g(x,x)Jx=ig(x,x)x,
\end{equation}
the eigenvalue of the eigenvector $x$ is $ig(x,x)$, which gives the result.
\end{dk}

\begin{lem}
There are two orbits of the adjoint action of $\Su(n,1)$ on rank 1 matrices in $\su(n,1)$. For a given
vector $x\in V$, the morphism $x\wedge Jx$ lies in one of the orbits, the morphism $-x\wedge Jx$ in the
other one.
\end{lem}
\begin{dk}
Any rank 1 morphism in $\su(n,1)$ has  all eigenvalues equal zero (if it would have an eigenvector with
non-zero eigenvalue, it would have to have at least one other, to be traceless; then it would not be
of rank 1). 
Then its canonical Jordan normal form has exactly one block of size 2.
According to the $\Su(n,1)$-orbits classification of $\su(n,1)$, see Lemma $\ref{orbity}$,
the morphism belongs to some of the type 2 orbit.

Denote $D=x\wedge Jx$. Then $Jx\wedge J(Jx)=D$ and $ax\wedge aJx = a^2J$, $a\in\R$, that is the morphism
$-x\wedge Jx$ cannot be written in the form $y\wedge Jy$ for any $y\in V$.

\end{dk}

\subsection{Geometric interpretation of $\wh C$ and $C$, the Lie algebra cone and its projectivization}

Let us follow the construction for $G=\Su(n,1)$.	
There is $\widehat{C}=Ad_{\Su(n,1)}(x)\subset
\su(n,1)$, where $x$ is some maximal root
element. 
The roots in $\su(n,1)$ have the same length 
and any rank 1 matrix in $\su(n,1)$ is a maximal root element.

\begin{lem}
The cone $C\subset\su(n,1)$ is isomorphic to $\C^{n}\setminus \{0,\dots,0\}$
\end{lem}

\begin{dk}
Since any maximal root element in $\su(n,1)$ is a rank 1 matrix, then the cone $\widehat C$ is equal to one of
the two orbits of the adjoint action on the rank 1 matrices (according to the orbit of the generating element). 
According to the lemmas $\ref{Kruhova akce}$ and $\ref{Bezestopy}$ we have
\begin{equation*}
\widehat C\cong \{x\in\C^{n+1}|g(x,x)=0\}/S^1. 
\end{equation*}
We describe the isomorphism explicitly. Let us choose the basis $e_1$,\dots,$e_{n+1}$ of $\C^{n+1}$ in which the metric $g$ 
has the standard form $(e^1)^2+\cdots +(e^n)^2-(e^{n+1})^2$.
We can choose a unique representative $\by$ in
each orbit of $S^1$ of the zero length vectors such that
$\by=(x_1,\dots,x_{n+1})=({\bf x},x_{n+1})\in
\C^{n+1}$ with $x_{n+1}\in\R^+$ with respect to the basis $(e_i)$. If we denote $\|\;\|$ the standard metric
$(e^1)^2+\cdots+(e^n)^2$ in $\C^n$ we can write the representative in the form $(\bx,\|\bx\|)$, where
$\bx$ is an arbitrary vector in $\C^n$. Hence $\widehat C$ is isomorphic to $\C^n\setminus \{0,\dots,0\}$.
\end{dk}

\begin{lem}
The projectivized cone $C\subset\P_o(\su(n,1))$ is isomorphic to $S^{2n-1}$. 
\end{lem}

\begin{dk}
Consider the surjective map $\pi:\wh
C\to\C^n\backslash(0,\dots,0)$,
$(\bx,\|\bx\|)\mapsto (\frac{x_1}{\|\bx\|},\dots, \frac{x_n}{\|\bx\|})$, which maps the elements in $\widehat C$ 
which differ by a real multiple to the same point on the sphere $S^{2n-1}$.
\end{dk}

\subsubsection{The module structure of the cone}

Being a homogeneous space (see \eqref{Homogeneous spaces}), the cone $C$ is a $SU(n,1)$, resp.
$\su(n,1)$-module. 

Let us notice, that the identification of the cone $C$ with the
sphere $S^{n+1}$ is subject to the choice of the standard hermitian
form on $\C^{n+1}$, respectively to the choice of the normal base with respect to it. We have to have
this in mind when considering different block forms of the matrices in $\su(n,1)$.

The proof of the lemma \eqref{orbity} shows that  if a matrix of the form
$\begin{pmatrix}
A&0\\
0&-\tr(A)
\end{pmatrix}$
is in $\su(n,1)$, then $A\in\un(n)$ (with respect to some orthonormal basis of the given hermitian form).

The one-parametric subgroup generated by $A$ are then matrices of the form

$
\begin{pmatrix}
G&0\\0&1\over \det(G)
\end{pmatrix}
$, where $G\in \Un(n)$. 
As we have seen above, a point in our sphere $S^{2n-1}$ corresponds to a class of null vectors in $\C^{n+1}$
which differ by a complex multiple. According to the lemma $\ref{invariance}$ the group $SU(n,1)$ acts on these
classes in a standard way and we get
\begin{eqnarray}
\label{Akce na sfere}
G\cdot\bx&\sim& \begin{pmatrix}G&0\\0&1\over \det(G)\end{pmatrix} (\bx,1)=\left(G\bx,\frac1{\det(G)}\right)\sim \det(G) G\bx
\end{eqnarray}
For the action of the matrices in $\un(n)$ of the above form on the tangent bundle of the sphere we get then:
\begin{equation}
\label{Akce algebry}
A\cdot \bx=(\tr(A)E+A)\bx,
\end{equation}
which is in accordance with the action of $\un(n-1)$ on $\C^{n-1}$ from $\eqref{akceun}$, if we represent
$A$ in the matrix form from $\eqref{matrixform}$.

Thus we get in fact the structure of $\Un(n)$ and $\un(n)$-module on $C$ regarded as a sphere $S^{2n-1}$.

The action of the whole group $SU(n,1)$ or the whole algebra $\su(n,1)$ respectively is then non-linear
(as a $\C^n\to\C^n$ mapping).

\begin{lem}
\label{Omega on C}
The canonical symplectic form on $\wh C\subset\g^*$ corresponds under our identification to the standard
symplectic form on $\C^n$, the one form $\lambda=\iota_{E_0}(\Omega)$ is then $\lambda=z_i\,d\overline{z_i}+
\overline{z_i}\,d z_i$ in the complex coordinates $z_i$ on $\C^n$.
\end{lem}
\begin{dk}
The Euler vector field on $\C^n$ is $E_0=z_i\,\frac{\del}{\del z_i}$. The form $\lambda=z_i\,d\overline{z_i}+
\overline{z_i}\,dz_i$ is then the only $U(n)$-invariant 
one-form $\alpha$ on $\C^n$ satisfying ${\cal L}_{E_0}\alpha=\alpha$ and $\alpha(E_0)=1$.
\end{dk}


Up to now, the construction was common for all the Bochner-Kaehler geometries. The choice of a transversal
symmetry of the canonical contact distibution actually determines the geometry.
Let $A\in \su(n,1)$ and let us consider the vector field $\xi{v}=\frac{\partial }{\partial t}|_0
\Ad(\exp(tA))v$ on $\widehat{C}$.
This vector field is a contact symmetry with respect to the distribution $\widehat{\D}$ on $\widehat{C}$
and thus it determines a section of $\widehat C\to C$ that is a contact form on $C$
(the identification of $\g$ and $\g^*$ gives an identification of $C\subset\su(n,1)$ and 
$C^*\subset \su(n,1)^*$).



The section $\lambda:C\to \widehat C$ is given by the equation $\lambda(\xi)=1$.
The image of $C$ in $\widehat C=\C^n-\{0\}$ is
then a hyperplane, which we will call $\Sigma_A$. The tangent space of $\Sigma_A$ is then characterized
by $T_v(C)=\{X\in \C^n|\w(X,iA\cdot v)=0\}$ and consequently there is
\begin{equation}
\label{Sigmadistribuce}
\D_v=\{X_v|g(X,iv)=0=g(X,iA\cdot v)\}.
\end{equation}

\odstavec{The projectivized cone $C$ and the $CR$-sphere} 
According to the lemma $\eqref{Homogeneous spaces}$, the projectivized cone $C$ is a homogeneous space $SU(n,1)/P$, 
where P is a parabolic subgroup of $SU(n,1)$, corresponding
to the subalgebra $\un(n)\oplus\C^{n-1}\oplus e_+^2$ of $\su(n,1)$, see $\eqref{Gradace algebry}$. As we have
seen in the previous lemma, it is isomorphic to the sphere $S^{2n-1}$. The adjoint action of $SU(n,1)$ on $\widehat C$ 
corresponds to the standard action of $SU(n,1)$ on the null-vectors (with respect to the standard hermitian form of 
the signature $(n,1)$) in $\C^{n+1}$, and thus as a homogeneous space it is exactly the $CR$-sphere 
(see for example \cite{Parabook}).

The underlying structure connected with this space is the canonic $CR$-distribution ${\cal D}$
on the sphere:

\begin{equation*}
\widehat{\D}_v=\{x\in \C^n|(x,iv)=0, (x,v)=0\}, 
\end{equation*}

where $(,)$ denotes the standard metric in $\C^n$.

\subsection{Classification of Bochner-Kaehler metrics}
 
All
Bochner-Kaehler manifolds come from the mentioned construction for the Lie
algebra $\su(n,1)$ and the resulting manifolds are isomorphic if we take in the
course of construction matrices $A_0$ lying on the same adjoint orbit of
$\Su(n,1)$ in $\su(n,1)$. Thus we can classify the Bochner-Kaehler manifolds
according to which orbit of the action induces the given manifold.

 There are four types of orbits of the adjoint action of $GL(n+1,\C)$ on
  $\su(n,1)$. Three types of these orbits are $SU(n,1)$ orbits as well, the fourth
one splits into two $SU(n,1)$ orbits.
   We describe the orbit types according to the Jordan blocks
   of the matrices in the orbits.

\begin{lem}
\label{orbity}
There are five types of orbits of the adjoint action of the $\Su(n,1)$ on $\su(n,1)$. If we represent a morphism
in the $\su(n,1)$ with a matrix $A$, than the orbit types look
as follows:

\begin{description}
\item[1.] The matrix is diagonizable  and its  eigenvalues are purely imaginary.
\item[2.] The eigenvalues of matrices in the orbit  are pure imaginary and there is
     just one Jordan block of the dimension 2 (there are $n$ eigenvectors). There exists
     an eigenvector $e$ and a root vector $f$, both in the block, such that 
      \begin{description}
\item[2a.] $(e,f)=i$.
\item[2b.] $(e,f)=-i$.
\end{description}
\item[3.] The eigenvalues of matrices in the orbit  are pure imaginary and there is
     just one Jordan block of the dimension 3 (there are $n-1$ eigenvectors).
\item[4.] There are $n-1$ pure imaginary eigenvalues corresponding to $n-1$
     eigenvectors and two eigenvalues $\lambda=\lambda_1+i\lambda_2$ and
     $\mu=\mu_1+i\mu_2$ 
\end{description}
\end{lem}
\begin{dk}
 It is easy to see that if $\lambda$ is an eigenvalue of a $\su(p,q)$-morphism then
  $-\overline\lambda$ is the eigenvalue of the morphism as well. 
There exists no nullplane in $V$, that is there are no two null-vectors $x$, $y$ in $V$ such that also
$h(x,y):=(x,y)=0$.
  
Let us suppose that the $\su(n,1)$-morphism $A$ has a Jordan block corresponding to an eigenvalue $\la$ with
non-zero real part, that is of size at least 2. Then there are $x$, $y\in V$ such that
$Ax=\la x$, $Ay=\la y+x$, $(x,x)=0$. Then
\begin{eqnarray*}
0=(Ax,y)+(x,Ay)=(\la x,y)+(x,\la y+x)=\la (x,y)+\overline\la(x,y)+(x,x)=2\Re\la(x,y),
\end{eqnarray*}
that is $(x,y)$=0. Further
\begin{equation*}
0=(Ay,y)+(y,Ay)=(\la y+x,y)+(y,\la y +x)=\la (y,y)+(x,y)+\overline\la(y,y)+(y,x)=2\Re\la(y,y),
\end{equation*}
which implies $(y,y)=0$ and we get a null-plane, which is a contradiction. All Jordan blocks corresponding
to the eigenvalues with non-zero real part have the size one.

Further if a $\su(n,1)$ morphism has an eigenvalue with non-zero real part, then it has to have at least
two eigenvalues which are not in $i\R$, otherwise the trace could not be zero. 
Let $\l=\la_1+i\la_2$ and $\mu=\mu_1+i\mu_2$, $\la_1$, $\mu_1$, $\la_2$, $\mu_2\in \R$, $\la_1\ne 0\ne \mu_1$ 
be eigenvalues of a $\su(n,1)$ morphism $A$, with non-zero real parts. Let $x$ resp. $y\in V$ be
the corresponding eigenvectors. Then
\begin{eqnarray*}
0=(Ax,x)+(x,Ax)=2(\la x,x)=2\la_1(x,x)
\end{eqnarray*} 
The eigenvector $x$ is then a null-vector. With the same argument is $y$ a null-vector too. 
\begin{eqnarray*}
0&=&(Ax,y)+(x,Ay)=(\la x,y)+(x,\mu y)=\la (x,y)+\overline\mu (x,y)\\
&=&(\la +\overline\mu)(x,y),
\end{eqnarray*}
and $\la=-\overline\mu$ otherwise $x$ and $y$ would generate a nullplane.

A $\su(n,1)$ morphism has consequently either exactly two eigenvalues $\la$, $\mu$ with non-zero real parts
and in this case $\la=-\overline \mu$, or all its eigenvalues are pure imaginary (or zero).
Let $V_{\la}$ and $V_{-\overline\la}$ be the eigenspaces corresponding to the eigenvalues with non-zero real parts.
The space $V_\la\oplus V_{-\overline\la}$ is $A$-invariant and so is the space $W:=(V_\la\oplus V_{-\overline\la})^\perp$:
let $x\in V_\la\oplus V_{-\overline\la}$ and $y\in W$  and then
\begin{equation}
(Ay,x)=-(y,Ax)=0.
\end{equation}
The restriction of the form $h$ to $V_\la\oplus V_{-\overline\la}$  has the signature $(1,1)$, therefore
is $h$ on $W$   positive definite, that means $A|_W\in \un(W,h)$ and $A|_W$ is diagonalizable (with
eigenvalues in $i\R$). The matrix form of the morphism $A$ is then $A=\diag(\la_1+i\la_2,-\la_1+i\la_2,i\la_3,
\dots,i\la_{n+1})$
in an apropriate basis
$e_1$,\dots $e_{n+1}$ such that 
\begin{equation*}
(e_1,e_1)=(e_2,e_2)=0,\quad (e_1,e_2)=1,\quad (e_1,e_i)=(e_2,e_i)=0\ \forall i\ge 3,\quad (e_i,e_j)=\delta_i^j\ 
\forall i,j\ge 3.
\end{equation*} 
Any such two morphisms lie evidently on the same $\Un(n,1)$ orbit, that is also on the same $\Su(n,1)$ orbit.

Let us check how big can be the Jordan blocks 
of the Jordan normal forms of the morphism $A$ with only pure imaginary (or zero) eigenvalues. 

\noindent
{\bf 1.} There exists  a Jordan block of of the size at least 4. Let  
$i\la$, $\la\in\R$ be on the diagonal . The
morphism $i\la$ is an unitary one, and so is the morphism $B:=A-i\la$. The morphism $B$ is then unitary and
in its Jordan normal form has the block of the form 
$\begin{pmatrix} 0&1&0&0\\0&0&1&0\\0&0&0&1\\0&0&0&0
\end{pmatrix}$. And let $e_1$, $e_2$, $e_3$, $e_4$ be the root vectors corresponding to the block, that is
$Be_4=e_3$, $Be_3=e_2$, $Be_2=e_1$, $Be_1=0$. Then we have
\begin{eqnarray*}
0&=&(Be_1,e_2)+(e_1,Be_2)=0+(e_1,e_1)\quad\impl (e_1,e_1)=0,\\
0&=&(Be_1,e_3)+(e_1,Be_3)=0+(e_1,e_2)\quad\impl (e_1,e_2)=0,\\
0&=&(Be_2,e_3)+(e_2,Be_3)=(e_1,e_3)+(e_2,e_2)\quad\impl (e_1,e_3)=-(e_2,e_2),\\
0&=&(Be_3,e_3)+(e_3,Be_3)=(e_2,e_3)+(e_3,e_2)=2\Re (e_2,e_3)\quad \impl (e_2,e_3)\in i\R \\
&&\mbox{(it will be used in the point 2)}\\
0&=&(Be_1,e_4)+(e_1,Be_4)=0+(e_1,e_3)\quad \impl (e_1,e_3)=0\quad\impl (e_2,e_2),
\end{eqnarray*}
and we would get a null-plane $(e_1,e_2)$. Thus there is no morphism in $\un(n,1)$ with the Jordan block of the size greater
then 3 (and so is no such morphism in $\su(n,1)$.

\noindent
{\bf 2.} Let there be a Jordan block of size 3. As in the first point we get vectors $e_1$, $e_2$, $e_3$ with
$(e_1,e_1)=(e_1,e_2)=0\stackrel{\mbox{\tiny sign. $(n,1)$}}{\impl} (e_2,e_2)>0$ and we may choose $e_2$ such
that $(e_2,e_2)=1$, then $(e_1,e_3)=-1$ and $(e_2,e_3)\in i\R$.

Now let us consider the transformation $f_1:=ae_1$, $f_2:=ae_2+be_1$, $f_3:=ae_3+be_2+e_1$, $a$, $b$, $c\in \C$.
Then we have still $Bf_1=0$, $Bf_2=f_1$, and $Bf_3=f_2$. After and easy computation one can choose $a$, $b$ and $c$ so that
\begin{equation*}
(f_1,f_1)=(f_3,f_3)=0,\quad (f_1,f_3)=-1,\quad (f_2,f_j)=\delta^j_2,\quad \mbox{for}\quad j=1,2,3
\end{equation*}
Now the space $\spa(e_1,e_2,e_3)$ is $A$-invariant  and the restriction
of $h$ to this space has signature $(2,1)$.
It follows that the space $W:=\spa(e_1,e_2,e_3)^\perp$ is $A$-invariant as well, and $h$ is positive definite
on it, that is $A|_W$ is diagonizable with unitary transformation of $W$ and the eigenvalues lie in $i\R$.

Again, all such $\su(n,1)$ morphisms lie on the same $\Su(n,1)$ orbit. 

\noindent
{\bf 3.} There exists a Jordan block of size 2. Let us suppose first, that there are at least two
blocks of size 2. Let $e_1$, $e_2$, $f_1$, $f_2$ be corresponding linearly independent vectors such that
$Ae_j=i\la_je_j$ and $Af_j=i\la_jf_j+e_j$, $j=1,2$. Then
\begin{equation*}
0=(Ae_j,f_j)+(e_j,Af_j)=i\la_j(e_j,f_j)+\overline{i\la_j}(e_j,f_j)_(e_j,e_j)=(e_j,e_j), \quad j=1,2.
\end{equation*}
Thus $(e_1,e_2)\ne 0$ according to the nullplane argument, and we can write
\begin{equation*}
0=(Ae_1,e_2)+(e_1,Ae_2)=i\la_1(e_1,e_2)+\overline{i\la_2}(e_1,e_2)=i(\la_1-\la_2)(e_1,e_2),
\end{equation*}
which implies $\la_1=\la_2=:\la$. Then
\begin{equation*}
0=(Ae_1,f_2)+(e_1,Af_2)=i\la(e_1,f_2)+\overline{i\la}(e_1,f_2)+(e_1,e_2),
\end{equation*}
and we get the contradiction. Thus there is just one block of size 2. Let $e_1$, $f_1$, $e_3$,\dots,$e_{n+1}$
be a basis such that $Af_1=i\la_1f_1+e_1$, $Ae_j=i\la_je_j$, $\la_j\in \R$ for all $j=1,\dots n+1$. Then there is
\begin{eqnarray}
0&=&(Ae_1,e_j)+(e_1,Ae_j)=i\la_1(f_1,e_j)+(e_1,e_j)+\overline{i\la_j}(f_1,e_j)\\
\label{prvnirovnice}
&=&i(\la_1-\la_j)(f_1,e_j)+(e_1,e_j),\\
\label{druharovnice}
0&=&(Ae_1,e_j)+(e_1,Ae_j)=i\la_1(e_1,e_j)+\overline{i\la_j}(e_1,e_j)=i(\la_1-\la_j)(e_1,e_j).
\end{eqnarray}
Then $\eqref{prvnirovnice}$ and $\eqref{druharovnice}$ imply according to the excluded third principle 
$(e_1,e_j)=0$ for $j=1,3,$\dots,$n+1$. Since $h$ is non-degenerate we have $(e_1,f_1)\ne 0$ and
\begin{equation*}
0=(Af_1,f_1)+(f_1,Af_1)=2\Re(Af_1,f_1)=2\Re(i\la_1(f_1,f_1)+(e_1,f_1))=2\Re(e_1,f_1),
\end{equation*}
that means $(e_1,f_1)\in i\R\setminus 0$. Without loss of generality we can suppose that $(e_1,f_1)=\epsilon i$,
$\epsilon=\pm 1$. If we consider the transformation $e_2:=f_1+ce_1$, then $Ae_2=i\la_1e_2+e_1$ and we
can choose $c\in\C$ so that $(e_2,e_2)=1$. 
Then again as in the previous points we can choose a basis $e_1$,\dots, $e_{n+1}$ ($e_1$ and $e_2$ are already
given) such that
the morphism $A$ is in its canonical Jordan form with exactly one Jordan block of size 2 in the basis,
 and for the basis vectors there is:
\begin{eqnarray*}
(e_1,e_1)&=&(e_2,e_2)=0, (e_1,e_2)=\epsilon i,\epsilon\in\{\pm 1\}, (e_1,e_j)=(e_2,e_j)=0,\quad \mbox{for}\quad j\ge 3\\
(e_j,e_k)&=&\delta_j^k\quad \mbox{for}\quad j,k\ge 3.
\end{eqnarray*} 
There are two orbits of the $\Su(n,1)$ of the morphisms of this type according to $\epsilon$, that is according
to the scalar product of $e_1$ and $e_2$ in the canonical basis. The number $\epsilon$ is evidently a $\Su(n,1)$
invariant.

{\bf 4.} The morphism $A$ is diagonalizable. Again, we can choose a basis $e_1$,\dots $e_{n+1}$ in which
the morphism has the diagonal form and
\begin{equation*}
(e_1,e_1)=-1,\quad (e_1,e_j)=0, \forall j\ge 2,\quad (e_j,e_k)=\delta^k_j\quad \forall j,k\ge 2,
\end{equation*}
and all such morphism lie on one $\Su(n,1)$ orbit.

In low dimensional cases ($\dim V\le 2$) we have then less orbits types. For $\dim V=3$ the morphism with
the Jordan block of size 3 has all the eigenvalues zero.
\end{dk}

\odstavec{Characteristic polynomial of the metrics} 
The characteristic polynomial determines all invariants of the adjoint orbit of the matrix $A\in \g$\footnote{
Let $\g\subset \End(V)$ be an irreducible representation of the Lie algebra $g$. Let
$\phi:\g\to \R(t)$, such that $\phi(\Ad_gx)=\phi(x)$ for all $x, g\in\g$ and $\gr(\phi(x))\le n$
Then $\phi(x)$ is a constant multiple of the characteric polynomial of $x$.} Thus the different types of adjoint orbits 
correspond to different
types of characteristic polynomials (distinguished according to their roots) and we get invariants of
the equivalent classes of Bochner-Kaehler metric.

The characteristic polynomial $p_A$ of the matices in $\Gamma_A$ is  according to $\eqref{matrixform}$:
\begin{eqnarray}
&\det
\begin{pmatrix}
\rho-\frac{1}{n+2}(\tr\rho)\I_n&u&u\\
-u^*&-\frac12(\tr\rho-i(f+1))&\frac i2(f-1)\\
u^*&\frac i2(1-f)& -\frac12(\tr\rho+i(f+1))
\end{pmatrix}- t\I_{n+2}\\
=&\det(\rho-\frac{1}{n+2}(\tr\rho)-t\I_n)(t^2+(\tr\rho)t+f+\frac14(\tr\rho)^2)+u^*\Cof(\rho-\frac{1}{n+2}(\tr\rho)-t\I_n)u,
\end{eqnarray} 
where $\Cof(X)$ means the cofactor matrix of $X$.

This is in accordance with the Bryant's result (see \cite{Bochner}) 
\footnote{The functions $\rho$, $u$, $f$ from \ref{thm:canonicalconn} correspond to Bryant's functions
$S$, $T$, $U$ as follows: $\rho=iS$, $u=-T$, $U=-f$.}
on the orbits of the diagonalizable
matrices with pure imaginary eigenvalues.

\subsection{One "nice" type of Bochner-Kaehler metrics}

We are now going to describe the first one of the five mentioned types of Bochner-Kaehler metrics in more detail.
Namely let us investigate those metrics which come from the construction if we take in the course of it
the matrix generating transversal symmetry to be diagonalizable with all eigenvalues pure imaginary (or zero). 
These are the matrices which acts as linear morphism ($U(n)$-morhpism actually) on the sphere $S^{2n-1}$
regarded as a projectivized cone $C$.

Any diagonalizable matrix  $\wh A\in\su(n,1)$ can be written in the form
$\diag(i\lambda_1,\dots,i\lambda_n,-i\sum_{i=1}^n\lambda_i)$,
where
$A=\diag(i\lambda_1,\dots,i\lambda_n)$ is in $\un(n)$. 

Then according to $\eqref{Akce algebry}$, $A$ acts on the tangent bundle of the spere:

\begin{equation*}
A\cdot
(x_1,\dots,x_n)=(i(\lambda_1+\sigma)x_1,\dots,i(\lambda_n+\sigma)x_n)\in
TC, \quad\sigma=\sum_{i=1}^n\lambda_i. 
\end{equation*}
The action of a matrix $A$ on $C^n$
thus
corresponds
to the multiplication with the matrix $A'$, where $A'=A+\tr(A)I$, $I$
being the indentity matrix.

Let $A_0\in\su(n,1)$, 
$A_0=
\begin{pmatrix}
A_0&0\\0&-\tr(A_0)
\end{pmatrix}$ 
be of the above diagonal form such that $\xi_0{v}=A_0\cdot v$ defines a transversal symmetry
on $C_0$, a non-empty open subset of $C$, that is $\xi_0{v}\notin \D_v$ on $C_0$.
This symmetry then defines a section $\lambda$ of the line bundle $\wh C\to C$. Then
$\lambda(C)=\Sigma_A\subset \wh C$, 
\begin{eqnarray}
\label{Sigmaeq}
\Sigma_A&=&\{(x\in \C^n|\lambda(\xi_0(x))=1\}\nonumber\\
&=&\{(x\in \C^n|(\xi_0(x),Jx)=1\}\nonumber\\
&=&\{(x\in \C^n| (A_0\cdot x, Jx )=1\}\nonumber\\
&=&\{x\in \C^n| (x,-(A_0+\tr(A_0))Jx)=1\}\nonumber\\
&=&\{x\in \C^n| \sum_{i=1}^n (\lambda_i+\sigma)|x_i|^2=1\}.\nonumber\\
\end{eqnarray}


Since the globalization brings to life some new demanding questions, we keep on  working 
locally only. Consider $U$ a regular open subset of $C$ with respect to $\xi_0$, that is
there
is a submersion $\pi_U:U\to M_U$ onto some manifold $M_U$, the set of leaves of
the foliation generated on $U$ by $\xi$. The whole of $C$ can be covered by
regular subsets. We write $M_U=U/T$.

Our goal is to determine the Bochner-Kaehler connection on $M_U$ which is induced
there acording to the general construction of special symplectic geometries.
The Bochner-Kaehler connections are with one-to-one correspondence with
the Bochner-Kaehler metrics which are further in one-to-one corresponce with
the pair consisting of the fundametal form of the Kaehler structure and the complex
structure on $M_U$.

There is the unique symplectic (which turns to the fundamental Kaehler one with the complex structure
on $M_U$) form $\w_U$ on $M_U$ such that
the pull-back of this form to $U\times \R\subset \wh C\subset \C^n$ is the
canonical symplectic form on $\C^n$ (this is the form that comes with the above
identifications from the Cartan-Killing form on $\su(n,1)$).

\begin{lem}
The complex projective space $\C P^n$ comes from our construction 
for $\g=\su(n+1,1)$ and
$$A=\left(\begin{array}{cccc}
-\frac{i}{2(n+2)}&&&\\
&\ddots&&\\
&&-\frac{i}{2(n+2)}&\\
&&&\frac{i(n+1)}{2(n+2)}
\end{array}\right).$$ 
\end{lem}

\begin{dk}
As we have already mentioned, the canonical $CR$-distribution on $S$ is the structure
which determines 
The matrix A acts on the sphere $S^{2n-1}$ according to \eqref{Akce algebry}
as linear map given by the matrix 
$$\left(\begin{array}{ccc}
-\frac{i}{2}&&\\
&\ddots&\\
&&-\frac{i}{2}\\
\end{array}\right),$$ 
and thus the contact symmetry $\xi$ is given as $\xi(p)=-\frac12ip$.
The hyperplane $\Sigma_A$ is a ${2n+1}$-dimensional sphere in $C^{n+1}$, so that it coincides
with the projectivization $C$ of the Lie algebra cone and the Kaehler form is the standard
form on $\C^n$. We get actually the cone construction from the theorem $\ref{Cone construction}$.

If we consider the matrix $A$ in the form of $\ref{matrixform}$, we see that the structure functions 
are $\rho=-\frac12J$, $f=1$, $u=0$.
\end{dk}

\odstavec{The generating metric on the contact distribution}
The Bochner-Kaehler structure on the sphere $S^{2n-1}$ lifts to the structure on $\Sigma_A$. 
The canonical complex structure on the contact distribution on $C$, that is canonical $CR$-distribution
on the sphere lifts to the complex structure on the distribution $\D_\Sigma\subset T\Sigma_A$, that
is to the complex structure on the contact distribution on the section $\Sigma_A$ of $\widehat C\to C$:
\begin{eqnarray}
  J_M(X)=JX-(X,A_0p)p-\frac{(X,p)}{|p|^2}Jp,
\end{eqnarray}
for vectors $X$, $Y\in T_p\Sigma$, and $|p|^2=(p,p)$.
This gives us then the metric on $\D_\Sigma$ 
\begin{eqnarray}
  \label{Sigmametric}
  g(X,Y)=\w(J_MX,Y)=(J_MX,JY)=(X,Y)-\frac{(X,p)(Y,p)}{|p|^2}.
\end{eqnarray}

This metric factors to the metric on $M_U\cong T\Sigma$, which is, 
according to
the construction, Bochner-Kaehler. This can be
confirmed also with the direct computation. 
We can view the metric $\eqref{Sigmametric}$ 
as a degenerated metric on the whole $\Sigma$. The corresponding Levi-Civita
connection is then
\begin{eqnarray}
  \label{Sigmaconnection}
  (\nabla_XY)_p=\nabla^0_XY+g(X,Y)\eta+(X,A_0JY)p,
\end{eqnarray}
where $X,Y,Z\in \frak X(\Sigma_A)$ $\nabla^0$ is a flat connection in $\C^n$, 
$\eta=A_oJp-|A_0p|^2p$, $p\in\Sigma_A$. For the Levi-Civita connection of the metric on $M_U$ we
have then
\begin{eqnarray*}
\label{Sigmaconnection2}
    \ov{\nabla_XY}=\nabla_{\ov X}\ov
    Y-\frac{1}{2}\w(X,Y)\xi_0+\alpha(X)\ov{J_MY} +  \alpha(Y)\ov{J_MX},
\end{eqnarray*}
 
where $X,Y\in \frak X(M)$, $\ov X$ is a lift of a vector field $X$ on $\frak X(M)$ to a vector
field on $\Sigma_A$ ($\ov X\in\D_\Sigma$), and $\alpha(X):=g(\ov X,\xi_0)=(\ov X,\xi_0)$. 

\odstavec{The Bochner-Kaehler form of the curvature of the metric.}
Let us
define the mapping $\rho:TM_U\to TM_U$ as 
\begin{eqnarray}
\label{rhomap}
\ov{\rho X}:=\nabla_{\xi_0}\ov X=
\nabla_{\ov X}\xi_0=A_0\ov X+g(\ov X,\xi_0)\eta+g(\ov X,A_0^2Jp)p.
\end{eqnarray}

\begin{lem}
The maping $\rho$ from \eqref{rhomap} is in $\un(n)$.
\end{lem}
\begin{dk}
Let us show first, that $\rho$ is well-defined, that is $g_0(\rho\oX,Jp)=0$:
\begin{eqnarray*}
(\overline{\rho X})&=&g(\rho X,Jp)=g(\nabla_{\xi_0}\oX,Jp)=\xi_0\underbrace{g(\oX,Jp)}_{=0}g(\oX,\nabla_{\xi_0}Jp)\\
&=&0-g(\oX,\nabla_{\xi_0}^0Jp+g(\xi_0,Jp)+g(\xi_0,-A_0p)p)\\
&=&-g(\oX,A_0Jp+(A_0p,Jp)A_0Jp)\\
&=&-g(\oX,A_0Jp-\underbrace{(p,A_0Jp)}_{=1}A_0Jp)\\
&=&-g(\oX,A_0Jp-A_0Jp)=0
\end{eqnarray*}

Next we show $\rho J_M=J_M\rho$:
\begin{eqnarray*}
\overline{\rho J_MX}&=&A_0J_MX+(\overline{J_MX},A_0\rho)\eta+(\overline{J_MX},A^2_0Jp)p\\
&=&A_0(J\oX-(\oX,A_0p)p-\frac1{|p|^2}(\oX,p)Jp)+(J\oX-(\oX,A_0p)p-\frac1{|p|^2}(\oX,p)Jp,A_0p)\eta\\
&&+(J\oX-(\oX,A_0p)p-\frac1{|p|^2}(\oX,p)Jp,A_0^2Jp)p\\
&=&A_0J\oX-(\oX,A_0p)A_0p-\frac1{|p|^2}(\oX,p)A_0Jp-\underbrace{(\oX,A_0p)}_{=0}\eta
-(\oX,A_0p)\underbrace{(p,A_0p)}_{=0}\eta\\
&&+(J\oX,A_0^2J\oX)p+(\oX,A_0p)\underbrace{(A_0p,JA_0p)}_{=0}p
\frac1{|p|^2}(\oX,p)(A_0Jp,A_0Jp)p\\
&=&A_0J\oX-(\oX,A_0p)+\frac1{|p|^2}(\oX,p)\underbrace{(-A_0Jp+\eta+|A_0p|^2)}_{=0}+(\oX,A_0^2p)p\\
&=&A_0J\oX-(\oX,A_0p)A_0p+(\oX,A_0^2p)p.
\end{eqnarray*}
For the $J_M\rho$ we have
\begin{eqnarray*}
\overline{J_M\rho X}&=&J\overline{\rho X}-(\overline{\rho X},A_0p)p-\frac1{|p|^2}(\overline{\rho X},p)Jp\\
&=&J(A_0\oX+(\oX,\xi_0)\eta+(\oX,A_0^2Jp)p)-(A_0\oX+(\oX,\xi_0)\eta+(\oX,A_0^2Jp)p,A_0p)p\\
&&-\frac1{|p|^2}(A_0\oX+(\oX,\xi_0)\eta+(\oX,A_0^2Jp)p,p)Jp\\
&=&A_0J\oX+(\oX,\xi_0)J\eta+(\oX,A_0^2Jp)Jp+(\oX,A_0^2p)p\\
&&-(\oX,\xi_0)\left(\underbrace{(A_0Jp,A_0p)}_{=0}
-|A_0p|^2\underbrace{(p,A_0p)}_{=0}\right)
+\frac1{|p|^2}(\oX,A_0p)Jp\\
&&-\frac1{|p|^2}(\oX,A_0p)\left(\underbrace{(A_0Jp,p)}_{=1}-|A_0p|^2|p|^2\right)Jp
-\frac1{|p|^2}(\oX,A_0^2Jp)|p|^2Jp\\
&=&A_0J\oX+(\oX,A_0p)\left(-A_0p-|A_0p|^2Jp+ \frac1{|p|^2}\left(Jp-\frac1{|p|^2}(1-|A_0p|^2|p|^2)\right)Jp\right)\\
&&+(\oX,A_0^2p)p\\
&=&A_0J\oX-(\oX,A_0p)A_0p+(\oX,A_0^2p)p
\end{eqnarray*}
Finally we show $g(\rho X,Y) = (A_0\oX,\oY)$, which implies $g(\rho X,Y)=-g(\rho Y, X)$
\begin{eqnarray*}
g(\rho X,Y)&=& g(\overline{\rho X}, \oY)=g(A_0\oX +g(\oX,\xi_0)\eta+g(\oX,A_0^2Jp)p,\oY)\\
&=&g(A_0\oX+g(\oX,\xi_0)A_0Jp,\oY)\\
&=&(A_0\oX,\oY)-  \frac1{|p|^2}(A_0\oX,p)(\oY,p)+(\oX,\xi_0)\left(\underbrace{(A_0Jp,\oY)}_{=0}-
\frac1{|p|^2}\underbrace{(A_0Jp,p)}_{=1}(\oY,p)\right)\\
&=&(A_0\oX,\oY)+\frac1{|p|^2}(\oX,A_0p)(\oY,p)-\frac1{|p|^2}(\oX,A_0p)(\oY,p)\\
&=&(A_0\oX,\oY)
\end{eqnarray*}

\end{dk}

With the help of the map $\rho$, we can express
the curvature of the Levi-Civita connection  on $M_U$ in the Bochner-Kaehler
form $\eqref{krivost2}$.

\begin{prop}
Let $M_U$ be a Bochner-Kaehler manifold which comes from the general construction
in $\cite{CS}$ for $\h=\un(n-1)$, $\g=\su(n,1)$, and $A_0\in\su(n,1)$. Then the curvature $R$ of the
Bochner-Kaehler metric on $M_U$ is given by

\begin{eqnarray*}
  R=R_{\frac{1}{2}\rho+\frac{1}{4}|A_0'p|^2J}
\end{eqnarray*}
\end{prop}
\begin{dk}
Let $\nabla$, $R$ be Levi-Civita connection and curvature of the Bochner Kaehler metric induced on $M_U$, and
$\on$, $\oR$ the ones of the metric $g$ on $\D_{\Sigma}$, $\al(X)=g(\oX,\xi_0)$ as before.
There is
\begin{eqnarray*}
(\nabla_X\al)(Y)&=&X\al(Y)-\al(\nabla_XY)=\oX g(\oY,\xi_0)-g(\overline{\nabla_XY},\xi_0)\\
&=&g(\nabla)_{\oX}\oY,\xi_0)+g(\oY,\nabla_{\oX}\xi_0)-g(\nabla_{\oX}\oY-\frac12\om(X,Y)\xi_0+\al(X)\overline{J_MY}
+\al(Y)\overline{J_MX},\xi_0)\\
&=&g(\oY,\overline{\rho X})+\frac12\om(X,Y)g(\xi_0,\xi_0)-\al(X)g(\overline{J_MY},\xi_0)-\al(Y)g(\overline{J_MX},\xi_0)\\
&=&g(\overline{\rho X},\oY)-g(X,J_MY))|A_0p|^2-\al(X)\al(J_MY)-\al(Y)\al(J_MX)
\end{eqnarray*}
Further
\begin{eqnarray*}
d\al(X,Y)&=&(\nabla_X)\al)(Y)-(\nabla_Y\al)(X)\\
&=&g(\overline{\rho X},\oY)-g(\overline{\rho Y},\oX)-g(\oX,\overline{J_MY}|A_0p|^2+
g(\oY,\overline{J_MX})|A_0p|^2\\
&&-\al(X)\al(J_MY)+\al(Y)\al(J_MX)-\al(Y)\al(J_MX)+\al(X)\al(J_MY)\\
&=&-2g(\oX,\overline{\rho Y})-2g(\oX,\overline{J_MY})|A_0p|^2.
\end{eqnarray*}
For the curvature tensor of the metric $g$ we get then
\begin{eqnarray*}
\overline{R(X,Y)Z}&=&R(\oX,\oY)\oZ+\left(-2g(\oX,\overline{\rho Y})-2g(\oX,\overline{J_MY})|A_0p|^2\right)\overline{J_MZ}\\
&&+\left(\al(Y)g(\oX,\oZ)-\al(X)g(\oY,\oZ)\right)\xi_0+g(\oY,\overline{J_MZ}\overline{\rho X}-g(\oX,\overline{J_MZ})\overline{\rho Y}\\
&&+\al(J_MZ)\left(\al(X)\overline{J_MY}-\al(Y)\overline{J_MX}\right)+\al(Z)\left(\al(X)\oY-\al(Y)\oX\right)\\
&&+\al(Z)\left(\al(J_MX)\overline{J_MY}-\al(J_MY)\overline{J_MX}\right)\\
&&+\left(-g(\oX,\overline{\rho Z})-g(\oX,\overline{J_MZ})|A_0p|^2-\al(X)\al(J_MZ)-\al(Z)\al(J_MX)\right)\overline{J_MY}\\
&&-\left(-g(\oY,\overline{\rho Z})-g(\oY,\overline{J_MZ})|A_0p|^2-\al(Y)\al(J_MZ)-\al(Z)\al(J_MY)\right)\overline{J_MX}\\
&&-2g(\oX,\overline{J_MY})\overline{\rho Z}\\
&=&R(\oX,\oY)\oZ-2\left(g(\oX,\rho\oY)+g(\oX,\overline{J_MY})|A_0p|^2\right)\overline{J_MZ}\\
&&+\left(\al(Y)g(\oX,\oZ)-\al(X)g(\oY,\oZ)\right)\xi_0-g(\overline{J_MY,\oZ}\overline{\rho X}+g(\overline{J_MX,\oZ}\overline{\rho Y}\\
&&+\al(Z)\left(\al(X)\oY-\al(Y)\oX\right)\\
&&+\left(g(\overline{\rho X},\oZ)+g(\ol{J_MX},\oZ)|A_0p|^2\right)\ol{J_MY}-\left(g(\ol{\rho Y},\oZ)+
g(\ol{J_MY},\oZ)|A_0p|^2\right)\ol{J_MX}\\
&&-2g(\oX,\ol{J_MY})\ol{\rho Z}\\
\end{eqnarray*}
\begin{eqnarray*}
&=&R(\oX,\oY)\oZ-2\left(g(\oX,\ol{\rho Y})+g(\oX,\ol{J_MY})|A_0p|^2\right)\ol{J_MZ}\\
&&+\left(\al(Y)g(\oX,\oZ)-\al(X)g(\oY,\oZ)\right)\xi_0+\left(\ol{\rho X}\wedge \ol{J_MY}\right)\oZ-
\left(\ol{\rho Y}\wedge \ol{J_MX}\right)\oZ\\
&&+\al(Z)\left(\al(X)\oY-\al(Y)\oX\right)+|A_0p|^2\left(\ol{J_MX}\wedge\ol{J_MY}\right)\oZ-2g(\oX,\ol{J_MY})\ol{\rho Z}\\
&=&\left(\oY\wedge A_0J\oX-\oX\wedge A_0J\oY+|A_0p|^2\oX\wedge\oY\right)\oZ\\
&&-(g(X,Z)Y-g(Y,Z)X,A_0^2p)p-2\left(g(\oX,\ol{\rho Y}+g(\oX,\ol{J_MY})|A_0p|^2\right)\ol{J_MZ}\\
&&+\left(\al(Y)g(\oX,\oZ)-\al(x)g(\oY,\oZ)\right)\xi_0+\left(\ol{\rho X}\wedge\ol{J_MY}-\ol{\rho Y}\wedge\ol{J_MX}
\right)\oZ\\
&&+\al(Z)\left(\al(X)\oY-\al(Y)\oX\right)+|A_0p|^2(\ol{J_MX}\wedge\ol{J_MY})\oZ-2g(\oX,\ol{J_MY})\ol{\rho Z}\\
&=&\left(\oY\wedge A_0J\oX-\oX\wedge A_0J\oY+|A_0p|^2\oX\wedge\oY\right)\oZ+(g(Y,Z)\oX-g(X,Z)\oY,A_0^2p)p\\
&&-2\left(g(\oX,\ol{\rho Y})+g(\oX,\ol{J_MY}|A_0p|^2\right)\ol{J_MZ}+\al(Y)(\oX\wedge\xi_0)\oZ-\al(X)(\oY\wedge\xi_0)\oZ\\
&&+\left(\ol{\rho X}\wedge\ol{J_MY}-\ol{\rho Y}\wedge\ol{J_MX}\right)\oZ+|A_0p|^2(\ol{J_MX}\wedge\ol{J_MY})\oZ-
2g(\oX,\ol{J_MY})\ol{\rho Z}\\
&=&\left(\oY\wedge A_0J\oX-\oX\wedge A_0J\oY+|A_0p|^2\oX\wedge\oY\right)\oZ+\left(g(Y,Z)\oX-g(X,Z)\oY,A_0^2p\right)p\\
&&-2\left(g(\oX,\ol{\rho Y})+g(\oX,\ol{J_MY})|A_0p|^2\right)\ol{J_MZ}+\al(Y)(\oX\wedge\xi_0)\oZ-\al(X)(\oY\wedge\xi_0)\oZ\\
&&+\left(\ol{\rho X}\wedge\ol{J_MY}-\ol{\rho Y}\wedge\ol{J_MX}\right)\oZ+|A_0p|^2(\ol{J_MX}\wedge\ol{J_MY})\oZ-
2g(\oX,\ol{J_MY})\ol{\rho Z}\\
&=&\left(\oY\wedge A_0J\oX-\oX\wedge A_0J\oY+|A_0p|^2\oX\wedge\oY\right)\oZ+\left(\oY\wedge (X,A_0^2p)p\right)\oZ-
\left(\oX\wedge(\oY,A_0^2p)p\right)\oZ\\
&&-2\left(g(\oX,\ol{\rho Y})+g(\oX,\ol{J_MY})|A_0p|^2\right)\ol{J_MZ}+\left(\al(Y)\oX\wedge\xi_0-\al(X)\oY\wedge\xi_0\right)\oZ\\
&&+\left(\ol{\rho X}\wedge\ol{J_MY}-\ol{\rho Y}\wedge\ol{J_MX}\right)\oZ+|A_0p|^2(\ol{J_MX}\wedge\ol{J_MY})\oZ\\
&&-2g(\oX,\ol{J_MY})\ol{\rho Z}\\
&=&\oY\wedge\left(A_0J\oX+(\oX,A_0^2p)p-\al(X)\xi_0\right)\oZ-\oX\wedge\left(A_0J\oY+(\oY,A_0^2p)p-\al(Y)\xi_0\right)\oZ\\
&&+|A_0p|^2\left((\oX\wedge\oY+\ol{J_MX}\wedge\ol{J_MY})\oZ-2g(\oX,\ol{J_MY})\ol{J_MZ}\right)-2g(\oX,\ol{\rho Y})\ol{J_MZ}\\
&&+\left(\ol{\rho X}\wedge\ol{J_MY}-\ol{\rho Y}\wedge\ol{J_MX}\right)\oZ-2g(\oX,\ol{J_MY})\ol{\rho Z}\\
&=&(\oY\wedge\ol{J_M\rho X})\oZ-(\oX\wedge\ol{J_M\rho Y})\oZ+|A_0p|^2\left((\oX\wedge\oY+\ol{J_MX}\wedge\ol{J_MY})\oZ-
2g(\oX,\ol{J_MY})\ol{J_MZ}\right)\\
&&-2g(\oX,\ol{\rho Y})\ol{J_MZ}+\left(\ol{\rho X}\wedge\ol{J_MY}-\ol{\rho Y}\wedge\ol{J_MX}\right)\oZ-2g(\oX,\ol{J_MY})\ol{\rho Z}\\
&=&\left(\oY\wedge\ol{J_M\rho X}-\oX\wedge\ol{J_M\rho Y}+\ol{\rho X}\wedge\ol{J_MY}-\ol{\rho Y}\wedge\ol{J_MX}-
2g(\oX,\ol{\rho Y})\ol{J_M}-2g(\oX,\ol{J_MY})\ol{\rho}\right)\oZ\\
&&+|A_0p|^2\left((\oX\wedge\oY)+\ol{J_MX}\wedge\ol{J_MY}-2g(\oX,\ol{J_MY})\ol{J_M}\right)\oZ,
\end{eqnarray*}
that is
\begin{eqnarray*}
R(X,Y)&=&Y\wedge (J_M\rho X)-X\wedge (J_M\rho Y)+\rho X\wedge J_MY-\rho Y\wedge J_MX-2g(X,\rho Y)J_M\\
&&-2g(X,J_MY)\rho
+|A_0\rho|^2\left(\oX\wedge Y+J_MX\wedge J_MY-2g(X,J_MY)J_M\right)
\end{eqnarray*}
\end{dk}



\begin{ex}

For the complex projective space we get the following characteristic polynomial $p(t)$ for the class of the Bochner-Kaehler metrics
with the constant holomorphic curvature equal 1:
$$p(t)=\left(t+\frac{i}{2(n+2)}\right)^{n+1}\left(t-\frac{i(n+1)}{2(n+2)}\right),$$
which corresponds to the Bryant's one (see section 4.1.1. in $\cite{Bochner}$).

\end{ex}

\section{Bochner tower}

There is a question if one can embedd (locally) the Bochner-Kaehler
manifold into a Bochner-Kaehler manifold of two (real) dimensions greater.

Following the geometric interpretation of the construction of Bochner-Kaehler
metrics, we would like to have some embedding $\su(n,1)\vloz \su(n+1,1)$, which
would induce the embedding of corresponding Bochner-Kaehler manifolds.

Consider the embedding

\begin{equation}
G\mapsto 
\begin{pmatrix}
1&0\\0&G
\end{pmatrix}
\end{equation}  
of the Lie group $SU(n,1)$ into $SU(n+1,1)$ and the corresponding embedding 

$A\mapsto 
\begin{pmatrix}
0&0\\0&A
\end{pmatrix}:=B$ of the Lie algebras.

This embedding yields the embedding of corresponding cones (and their projectivizations)
in the Lie algebras, we have $C_{\su(n,1)}\equiv S^{2n-1}\vloz S^{2n+1}\equiv C_{\su(n+1,1)}$.
This embedding is evidently $\su(n,1)$-equivariant ($A\in\su(n,1)$ acts on $C_{\su(n,1)}$ according to the action $\eqref{Akce algebry}$,
$B$ acts on the embedded cone according to $\eqref{Akce algebry}$ for $\su(n+1,1)$).

A matrix $A\in\su(n,1)$ acts on $\C$ the same way as the matrix
\begin{equation}
\lower .2cm\hbox{\mbox{\strednimat D}}_{\lambda_0}=
\left(\begin{array}{c|c}
\lambda_0&0\hfill\dots\hfill0\\
\hline
0&\lower .2cm\hbox{\vbox to 0pt{\hbox{\Huge $A$\raise1.1mm\hbox{-}\normalsize\raise.2cm\hbox{$\frac{\lambda_0}{n+1}$}\Huge$E$\normalsize$_n$}}}\\
\vdots&\\
0&\\
\end{array}\right)\\
\end{equation}
on the image of $C_{\su(n,1)}\subset C_{\su(n+1,1)}$ under the described embedding.

Consequently we get the theorem

\begin{thm}
For any $\lambda_0\in \R$ and any $A\in \su(n,1)$ the Bochner-Kaehler manifold corresponding to $A$
(that is $T/\Sigma_A$, see $\ref{thm:canonicalconn}$) can be embedded totally geodesicly into to the Bochner-Kaehler
manifold $T/\Sigma_{D_{\lambda_0}}$, where $D_{\lambda_0}$ is given above.
\end{thm}

\begin{dk}
The manifold $\Sigma_{D_{\lambda_0}}$ is given by the equations $\eqref{Sigmaeq}$. Then the vector
$(0,v)\in \C^{n+1}$, $v\in \C^n$ lies in $\Sigma_{D_{\lambda_0}}$ evidently iff the vector $v$ lies in $\Sigma_A$. 
The contact distribution ${\cal D}_{\Sigma_{D_{\lambda_0}}}$ on $\Sigma_{D_{\lambda_0}}$ is given
by the equation $\eqref{Sigmadistribuce}$ and apparently
${\cal D}_{\Sigma_{D_{\lambda_0}}}\cap (0,\C^n)=(0,{\cal D}_{\Sigma_A})$. 
Moreover the action $\eqref{Akce algebry}$ shows that for a $v\in \Sigma_A$ there is 
$(0,\xi_A(v))=\xi_{D_{\lambda_0}}$. Thus the factor manifold $\Sigma_A/T_{\xi_A}$ is embedded into
the manifold $\Sigma_{D_{\lambda_0}}/T_{D_{\lambda_0}}$.

The lift of the Levi-Civita connection of the Bochner-Kaehler to the contact distribution $D$ on
$\Sigma_{D_{\lambda_0}}$ preserves the set of vectors of the form $(0,v)$, $v\in C^n$ and thus it preserves
also the distribution $(0,\Sigma_A)$.
Thus the described embedding is a totally geodesic one.
\end{dk}

\begin{rem}
The previous statement can be reformulated as: Bochner-Kaehler manifold can be totally geodesic ebedded into the one-parametric class of
Bochner-Kaehler manifolds of the (complex) dimension one higher. 
\end{rem}

\begin{ex}
Let us compute what the given embedding yields for the complex projective space $\C P^{n-1}$.
As we have seen, the complex projective space comes from the construction for
$\g=\su(n,1)$ and 

$A=\left(\begin{array}{cccc}
-\frac{i}{2(n+1)}&&&\\
&\ddots&&\\
&&-\frac{i}{2(n+1)}&\\
&&&\frac{in}{2(n+1)}
\end{array}\right)$. 
Choosing $\lambda_0=-\frac{i}{2(n+2)}$ we get $\C P^{n-1}$ we get 
\begin{equation*}
D_0=\left(\begin{array}{cccc}
-\frac{i}{2(n+2)}&&&\\
&\ddots&&\\
&&-\frac{i}{2(n+2)}&\\
&&&\frac{i(n+1)}{2(n+2)}
\end{array}\right),
\end{equation*}
and the Bochner-Kaehler manifold corresponding to this matrix is $\C P^n$. 
\end{ex}

\section{Bochner-K\"ahler and Ricci-type connections duality}

In this section  we describe the duality between the manifolds with the Bochner-K\"ahler metrics of
type 1. (see \ref{orbity}) and Ricci flat connections. 

Recall the general construction from the section 3. So far we were interested in the case with $\g_1=\su(n,1)$.
If we consider the construction for the parabolic 2-gradable algebra$\g_2:=\sp(n,\R)$, we get a manifold with the 
connection of Ricci type. 
Recall the two standard embeddings of $\un(n+1)$, first into $\su(n+1,1)$ (that was described in the previous
section),
second into $\sp(n+1,\R)$.

\begin{thm}
Consider the action of the Lie algebras $\g_1=\su(n+1,1)$ and $\g_2=\sp(n+1,\R)$ on the projectivized cones
$\widehat C_1$, $\widehat C_2$. Then the following are equivalent
\begin{itemize}
\item[i)] For $a_i\in\g_i$ the actions of $T_{a_i}\subset G_i$ on $\widehat C_i$ are conjugate for i=1,2.
\item[ii)] $a_i\in\un(n+1)$, where $\un(n+1)\subset\g_i$ for $i=1,2$ via the two standard embeddings.
\end{itemize}
\end{thm}
\begin{dk}
We have already computed the action of $\un(n+1)$ on $C_1$ (see \eqref{Akce algebry}). 
Observe, that the diagonalizable matrices in $\su(n+1,1)$ are  the only matrices, which act
on $C_1$ (which is isomorphic to the sphere $S^{2n+1}\subset\C^{n+1}$) in the standard way 
(as on the vectors in $\C^n$).
As for the action
on $C_2$ we have to go quickly through the general construction (Section 3) for $\g=\sp(n+1)$.
We have $S^2(\R^{2n+2})\cong \sp(n+1)$ ($(x\circ y)(z)=\om(x,z)y+\om(y,z)x$, for $x$, $y$, $z\in\R^{2n+2}$).
You can prove with an easy computation  as in the \ref{invariance}, that the map $x\mapsto x^2$ is
the $Sp(n+1)$-module homomorphism of the space $\R^{2n+2}/\Z_2$ and $S^2(\R^{2n+2})$,
where the image of the morphism corresponds to rank-one elements in $\sp(n+1)$.
The action of $\sp(n+1)$ on $\hat C_2\cong\R^{2n+2}/\Z_2$ is just a standard one, and  thus 
the action of $\un(n+1)\hookarrow\sp(n+1)$ on 
$C_2=\P_0(\hat C_2)\cong\R\P^{2n+2}$ is just a standard action on the real projective space. This is locally
the same as the action of $\un(n+1)\hookarrow\su(n+1,1)$ on the sphere $C_1$.
\end{dk}

Thus we get the following theorem:
\begin{thm}\hbox{}
\hfil\break

\vskip -5mm

\vbox{
\begin{itemize}
\item[i)] Let $(M,\om,\nabla)$ be a symplectic manifold with a connection of Ricci type, and suppose that the
corresponding element $A\in\sp(n+1,\R)$ from \ref{locallygenerating} is conjugate to an element of $\un(n+1)\subset
\sp(n+1,\R)$. Then $M$ carries a canonical Bochner-K\"ahler metric whose K\"ahler form is given by $\om$.
\item[ii)] Converselly, let $(M,J,\w)$ be a Bochner-K\"ahler metric such that the element $a\in\su(n+1,1)$
from \ref{locallygenerating} is conjugate to an element of $\un(n+1)\subset\su(n+1,1)$. Then $(M,\om)$ carries a canonical
connection of Ricci-type.
\end{itemize}
}
\end{thm}

\end{document}